\documentclass[review,onefignum,onetabnum]{siamonline190516}

\usepackage{amssymb}
\usepackage[linesnumbered,lined,boxruled,commentsnumbered]{algorithm2e}
\usepackage{bbm}
\usepackage{amsfonts}
\DeclareMathSymbol{\shortminus}{\mathbin}{AMSa}{"39}

\usepackage{lipsum}
\usepackage{amsfonts}
\usepackage{graphicx}
\usepackage{epstopdf}
\usepackage{caption}
\usepackage{subcaption}
\ifpdf
  \DeclareGraphicsExtensions{.eps,.pdf,.png,.jpg}
\else
  \DeclareGraphicsExtensions{.eps}
\fi

\usepackage{enumitem}
\setlist[enumerate]{leftmargin=.5in}
\setlist[itemize]{leftmargin=.5in}

\newsiamremark{remark}{Remark}
\newsiamremark{hypothesis}{Hypothesis}
\crefname{hypothesis}{Hypothesis}{Hypotheses}
\newsiamthm{claim}{Claim}

\headers{The self-consistent field iteration for $p$-spectral clustering}{P. Upadhyaya, E. Jarlebring, F. Tudisco}

\title{The self-consistent field iteration for $p$-spectral clustering}

\author{
Parikshit Upadhyaya\thanks{KTH Royal Institute of Technology, Stockholm, Sweden (\email{pup@kth.se})}
\and Elias Jarlebring\thanks{KTH Royal Institute of Technology, Stockholm, Sweden (\email{eliasj@kth.se}).}
\and Francesco Tudisco\thanks{GSSI Gran Sasso Science Institute, L'Aquila, Italy
  (\email{francesco.tudisco@gssi.it}).}
}

\usepackage{amsopn}
\DeclareMathOperator{\diag}{diag}
\DeclareMathOperator{\sign}{sgn}

\newcommand{\norm}[1]{\left\lVert#1\right\rVert}
\newcommand{\dd}{\mathcal{D}}
\newcommand{\rcut}{\operatorname{RCut}}
\newcommand{\cut}{\operatorname{Cut}}
\newcommand{\ncut}{\operatorname{NCut}}
\newcommand{\vol}{\operatorname{vol}}
\newcommand{\rcc}{\operatorname{RCC}}
\newcommand{\ncc}{\operatorname{NCC}}
\newcommand{\sm}{\shortminus}

\usepackage[final]{changes} 
\definechangesauthor[color=magenta]{FT}

\begin{document}

\maketitle

\begin{abstract}
The self-consistent field (SCF) iteration, combined with its variants, is one of the most widely used algorithms in quantum chemistry. We propose a procedure to adapt the SCF iteration for the $p$-Laplacian eigenproblem, which is an important problem in the field of unsupervised learning. We formulate the $p$-Laplacian eigenproblem as a type of nonlinear eigenproblem with one eigenvector nonlinearity \deleted[id=FT]{(NEPv1)}, which then allows us to adapt the SCF iteration for its solution after the application of \replaced[id=FT]{suitable}{certain}  regularization techniques. The results of our numerical experiments confirm the viablity of our approach. 
\end{abstract}



\section{Introduction}
In data science and machine learning, \textit{clustering} refers to the problem of dividing a given unlabelled dataset into relevant groups or clusters, guided by \replaced[id=FT]{suitable metrics that quantify}{a metric that quantifies} the pairwise similarity of the points in the dataset \added[id=FT]{as well as the overall quality of the cluster assignment}. Being one of the major approaches in unsupervised learning, it has a diverse range of applications, for example in the analysis of social networks \cite{NINA:2007:SOCIAL}, grouping a large number of correlated outputs from search queries \cite{YU:2020:PECOS}, medical imaging \cite{Selvan2017}, image segmentation \cite{SHI:2000} and biological sequence clustering \cite{REMM:2001}, to name a few. Clustering techniques can differ widely based on the features of the dataset that are considered to be relevant, the metric that is chosen to measure the pairwise \added[id=FT]{similarities} \deleted[id=FT]{similarity based on the relevant features} and the criterion that is optimized to eventually find the clusters. A comprehensive review of the many different clustering techniques can be found in \cite{XU:2015}.

\added[id=FT]{The term \textit{``spectral clustering''} refers to the family of clustering methods which are based on the spectrum of an operator defined on a graph whose edges model data similarities. 
With generic point cloud data, edges are typically based on $k$-nearest-neighbors or
$\epsilon$-neighbors,
but spectral clustering approaches are also often used directly on relational data coming from, e.g., social
networks or web graphs~\cite{fortunato2016community,vonLuxburg:2007:SPECLUST,ng2002spectral}. The different types of spectral clustering algorithms differ mainly in the choice of the graph operator.   Typical examples are the Laplacian $L=D-W$ and its normalized version $L_n = I-D^{\shortminus\frac{1}{2}}WD^{\shortminus\frac{1}{2}}$, where $D$ and $W$ are the degree and adjacency matrices of the graph, respectively. Other choices include the modularity matrix \cite{fasino2014algebraic,newman2006finding}, the Bethe Hessian \cite{krzakala2013spectral} and the heat kernel Laplacian \cite{kloster2014heat}, for example. While being widely explored in various areas, it is known that clustering methods based on linear graph operators may perform poorly on real-world data. To overcome this issue, nonlinear graph operators have been introduced in recent years as a generalization of more standard linear graph mappings \cite{boyd2018simplified,Hein:2009:PLAPLACIAN,tudisco2018community}. Among them, the $p$-Laplacian received considerable attention due to its simple definition, its connection  with Cheeger (isoperimetric) constants and nodal domains on graphs and manifolds, as well as its remarkable clustering performance \cite{Hein:2009:PLAPLACIAN,keller2016general,szlam2010total,tudisco2018nodal}. 

 From a graph-theoretic point of view, the Cheeger inequality connects the spectrum of graph Laplacians and $p$-Laplacians to the problem of minimizing balanced graph cuts and provides one of the main justifications for the use these operators for spectral clustering.}
For example, the second smallest eigenvalue of the unnormalized Laplacian $L$ sets a lower bound for the \textit{ratio cut} \cite{HAGEN:1991} of a graph $G=(V,E)$, which, for a partition of the node set $V$ into two clusters $C$ and its complement $C^c:=V\setminus C$, is defined as
\begin{equation}\label{eq:rcut}
\rcut(C,C^c) = \left(\frac{1}{|C|}+\frac{1}{|C^c|}\right)\cut(C,C^c) = \left(\frac{1}{|C|}+\frac{1}{|C^c|}\right)\sum\limits_{i\in C, j\in C^c}w_{ij},
\end{equation}
where $\cut(C,C^c)$ is the graph cut denoting the sum of all \added[id=FT]{the weights of} the edges that connect a node in $C$ with one in $C^c$. Similarly, the smallest eigenvalue of the normalized graph Laplacian $L_n$ sets a lower bound for the \textit{normalized cut} \cite{SHI:2000}, defined as
\begin{equation}
\ncut(C,C^c) = \left(\frac{1}{\vol(C)}+\frac{1}{\vol(C^c)}\right)\cut(C,C^c) = \left(\frac{1}{\vol(C)}+\frac{1}{\vol(C^c)}\right)\sum\limits_{i\in C, j\in C^c}w_{ij},
\end{equation}
where
\[
\vol(C) = \sum\limits_{i\in C}d_i
\]
 denotes the volume of the set $C$ obtained by summing the degrees of all the nodes in $C$. Minimizing the ratio or normalized cuts entails finding clusters of relatively \textit{balanced size}, where the size of a cluster \added[id=FT]{$C$} refers to the number of nodes \replaced[id=FT]{in $C$}{contained in it} or \added[id=FT]{to} its volume, respectively. A different balancing \added[id=FT]{function} is imposed if we minimize the ratio Cheeger cut ($\rcc$) or the normalized Cheeger cut ($\ncc$), which are defined as
 \begin{equation}\label{eq:rcc}
\rcc(C,C^c) = \frac{\cut(C,C^c)}{\min\{|C|,|C^c|\}},\;\; \ncc(C,C^c) = \frac{\cut(C,C^c)}{\min\{\vol(C),\vol(C^c)\}}.
 \end{equation}
In \cite{AMGHIBECH:2003} and \cite{Hein:2009:PLAPLACIAN}, the authors prove inequalitues that provide sharp bounds relating the second eigenvalue of the nonlinear $p$-Laplacian and optimal values of Cheeger cuts. If $\lambda_{p}^{(2)}$ and $\lambda_{p,n}^{(2)}$ are the second smallest eigenvalues of the unnormalized and normalized $p$-Laplacian respectively, then Amghibech shows in \cite{AMGHIBECH:2003} that
\begin{equation}
2^{p-1}{\left(\frac{h_{\ncc}}{p}\right)}^p \leq \lambda_p^{(2)} \leq 2^{p-1}h_{\ncc},
\end{equation}
and B\"uhler and Hein show in \cite{Hein:2009:PLAPLACIAN} that
\begin{equation}
{\left(\frac{2}{d_{\max}}\right)}^{p-1}{\left(\frac{h_{\rcc}}{p}\right)}^p \leq \lambda_{p,n}^{(2)} \leq 2^{p-1}h_{\rcc},
\end{equation}
where $h_{\ncc}$ and $h_{\rcc}$ are the optimal values of the normalized Cheeger cut and the ratio Cheeger cut and $d_{\max}$ is the maximum degree of all the nodes in the graph. In the limit $p\to 1$, these eigenvalues converge to the exact optimal value of the respective Cheeger cut, \added[id=FT]{showing that the spectrum of the $p$-Laplacian for $p\approx 1$ has superior clustering performance than its linear counter part}. 

\added[id=FT]{While being theoretically better suited for clustering purposes, computing the spectrum of the $p$-Laplacian for $p\neq 2$ is, in general, severely more challenging than the linear case $p=2$.} This has motivated several \replaced[id=FT]{authors}{papers in the literature} to focus on algorithms to compute eigenpairs of the $p$-Laplacian \replaced[id=FT]{for general $p\geq 1$}{, which cannot be formulated as a linear eigenvalue problem except for $p = 2$.} An inverse power method for the 1-Laplacian was proposed in \cite{Hein:2010:IPM}, \added[id=FT]{later extended into the RatioDCA algorithm and its generalizations \cite{hein2011beyond,jost2014nonlinear,tudisco2018community}.}  A modified gradient descent method is used to find eigenvectors that satisfy an orthogonality condition in \cite{Luo:2010}. In \cite{PASADAKIS:2020}, the authors modify the orthogonality condition and propose a Riemannian optimization algorithm on the Grassmanian associated with the orthogonality condition. A common theme in all these approaches is that they formulate the $p$-Laplacian eigenproblem as a non-convex optimization problem.

In this work, we look at the $p$-Laplacian eigenproblem from the perspective of numerical linear algebra and formulate it as what is known in the literature as a \textit{nonlinear eigenproblem with eigenvector nonlinearity} (NEPv) \cite{CAI:2018}. The NEPv is ubiquitous in the field of quantum chemistry, arising from numerical discretizations and approximations of the Schr\"odinger equation that lead to a nonlinear problem \cite{Saad:2010:ELECSTRUCT}. In general, the nonlinear terms in the quantum chemistry problems depend on the matrix of eigenvectors as a whole, where each eigenvector corresponds to a possible state the quantum system can occupy. However, there are certain applications like solving the Gross-Pitaevskii equation \cite{Bao:2004:BOSEEINSTEIN} where one is interested in computing just the ground state. In such cases, the nonlinear terms are functions of a single vector and we call them \textit{nonlinear eigenproblems with \textbf{one} eigenvector nonlinearity} (NEPv1). The NEPv1 is hence a specific type of NEPv, and the $p$-Laplacian eigenproblem can be formulated as a NEPv1, as we show in Section~\ref{sec:nepvform}.

A widely used approach to solve the NEPv1 in quantum chemistry is the \textit{self-consistent field iteration} (SCF), combined with its variants \cite{Saad:2010:ELECSTRUCT}. The SCF iteration is a fixed point iteration that starts with an initial guess for the eigenvector and in each step solves a linear eigenvalue problem constructed using the iterate from the previous step.  It has been studied extensively in the past few decades resulting in several convergence characterizations \cite{DENSITY:2021,Bai:2020:OPTIMALSCF,Cances:2000:SCF}, schemes to accelerate its convergence \cite{WALKER:ANDERSON:2011,Pulay:1980:CONVERGENCE} and make it more robust \cite{SAUNDERS:LEVEL:1993,ZERNER:DAMPING:1979}. This makes it a natural candidate to be used as an algorithm to compute eigenpairs of the $p$-Laplacian.

In this paper, we present two contributions to the field of numerical algorithms for the $p$-Laplacian:
\begin{itemize}
\item[(i)] Two different ways to formulate the $p$-Laplacian eigenproblem as a NEPv1.

\item[(ii)]  An adaptation of the SCF iteration to compute the eigenpair corresponding to the second smallest eigenvalue of the $p$-Laplacian, that is the eigenpair corresponding to the optimal Cheeger cuts.
\end{itemize}

Contribution (ii) consists of appropriately handling singularities of the problem. The situation $p\rightarrow 1$ is of main interest since we then have an
accurate approximation of the discrete optimization problem. However, as $p\rightarrow 1$ the nonlinear operator applied to a vector in general has elements of very large or very small magnitude. This is handled by replacing certain terms in the definition with smooth counterparts.

The rest of the paper is organized into four sections. In Section~\ref{sec:probform}, we provide a formal introduction to the $p$-Laplacian eigenproblem and show how it can be formulated in two different ways as a NEPv1. Section~\ref{sec:scfplap} introduces our numerical approach which is based on adapting the SCF iteration to a regularized form of a NEPv1-formulation of Section~\ref{sec:probform}. Section~\ref{sec:numex} contains the results of numerical experiments done by applying the technqiues introduced in  Section~\ref{sec:scfplap} to some standard test problems in machine learning. We end with some concluding remarks in Section~\ref{sec:conc}.

\section{Problem formulation}\label{sec:probform}
We assume that the dataset we wish to apply the clustering algorithm to is a collection of $n$ data points which can be denoted with indices in $\mathcal{I}_n = \{1,2,\ldots,n\}$. Additionally, we also assume that there is a notion of similarity $s:\mathcal{I}_n\times\mathcal{I}_n\to \mathbb{R}_{+}$, such that a higher value of $s(i,j)$ implies more similarity between the items $i$ and $j$. The existence of such a function $s$ allows us to construct an undirected weighted graph where the data points are the nodes. The weight of an edge connecting $i$ and $j$ depends on $s(i,j)$ but is not necessarily equal to it. A common choice for $s$ is the Gaussian kernel given by $s(i,j) = e^{-\norm{x_i-x_j}^2/\sigma^2}$, where $x_i$ and $x_j$ denote the coordinates of nodes $i$ and $j$, and $\sigma$ controls the width of the kernel.

Hence, we consider an undirected weighted graph $G:(V,E)$ with $n$ nodes in $V = \mathcal{I}_n$ and $m$ edges in $E$. The symmetric adjacency matrix is denoted by $W$, such that $w_{ij} \geq 0$ is the weight of an edge connecting node $i$ with node $j$.

We denote the degree matrix with $D$, which is a diagonal matrix such that $d_{i,i}$ contains the degree of node $i$, that is
\begin{equation}
d_{i,i} = \sum\limits_{j=1}^n w_{ij}.
\end{equation}

We assume the lexical order of the nodes is based on their indices in $V$. Using this lexical order of the nodes, we can assign a lexical order to the edges in $E$. We assign an edge $(i,j)$ a lower index than an edge $(k,l)$ if \textit{one of the following two conditions} is true:
\begin{enumerate}
\item $\min(i,j) < \min(k,l)$.
\item $\min(i,j) = \min(k,l)$ \textit{and} $\max(i,j) < \max(k,l)$.
\end{enumerate}

Based on this lexical order of the edges, we can associate each unordered pair in $E$ with an index in $\mathcal{I}_m$. For each edge $k$, we denote the index of the lesser (based on the order of the nodes) of the two nodes connected by the edge as $k_1$, and the other (the node with greater index) $k_2$. This lexical relationship between the edges and the corresponding nodes can be encoded in the graph \textit{incidence matrix} $B \in \mathbb{R}^{m\times n}$, defined as
\begin{equation}\label{eq:bdef}
B_{ki} = \begin{cases}
      -1,&\;\;\textrm{if}\;i = k_1,\\
      1,&\;\;\textrm{if}\;i = k_2,\\
      0,&\;\;\textrm{otherwise}.
    \end{cases}
\end{equation}
Each row in $B$ corresponds to an edge and contains zeros everywhere, except for the indices which correspond to nodes connected by the edge. We define $D^w \in \mathbb{R}^{m\times m}$ as a diagonal matrix with the $k$th diagonal entry denoting the weight of edge $k$, that is
\begin{equation}\label{eq:dw}
D^w_{k,k} = w_{k_1,k_2}.
\end{equation}

We define $\sign :\mathbb{R}^l\to\{-1,0,1\}^l$ to be the sign function, which applies elementwise to a vector.
For any vector $y\in\mathbb{R}^l$, $|y|$ denotes the vector containing absolute values of the elements in $y$.
\subsection{The $p$-Laplacian eigenproblem}\label{sec:plapeig}
The unnormalized graph Laplacian $L$ is defined as the difference of the degree and adjacency matrices, that is
\begin{equation}
L = D-W.
\end{equation}
Note that $L$ is a symmetric positive definite linear operator which, for any $x\in\mathbb{R}^n$ with each entry corresponding to a node, induces an inner product given by
\begin{equation}
x^TLx  = x^TDx - x^TWx = \sum\limits_{i,j\in V}w_{ij}(x_i-x_j)^2.
\end{equation}
 This is a quadratic function of the differences of elements in $x$ corresponding to the nodes that are connected by the respective edges. If we consider an inner product that is a function of these differences raised to $p$-th powers instead, the corresponding operator turns out to be a nonlinear generalization of $L$. This is how the $p$-Laplacian is defined implicitly.
\begin{definition}[$p$-Laplacian defined implicitly]
The unnormalized $p$-Laplacian is the operator $\Delta_p:\mathbb{R}^n\to\mathbb{R}^n$ which induces the inner product
\begin{equation}\label{eq:plapinn}
x^T\Delta_p(x) = \sum\limits_{i,j\in V}w_{ij}|x_i-x_j|^p.
\end{equation}
\end{definition}
Hence, we can refer to $L$ as the unnormalized 2-Laplacian. For the rest of the paper, whenever we refer to the $p$-Laplacian, we will mean the unnormalized $p$-Laplacian, unless specified otherwise.

In \cite{AMGHIBECH:2003}, it is shown that the explicit action of the $p$-Laplacian on $x$ produces a vector whose elements are given by
\begin{equation}\label{eq:plapaction}
(\Delta_p(x))_i =  \sum\limits_{j\in V}w_{ij}|x_i-x_j|^{p-1}\sign(x_i-x_j),\;\; i\in\mathbb{I}_n
\end{equation}
We define $\Phi_p:\mathbb{R}^l\to\mathbb{R}^l$, where $l$ can be different from $n$.
\begin{equation}\label{eq:phipdef2}
(\Phi_p(x))_i = |x_i|^{p-1}\sign(x_i),\;\; i = 1,\ldots, l.
\end{equation}
The action of the $p$-Laplacian can then be written in terms of $\Phi_p$ as
\begin{equation}\label{eq:plapexpli}
(\Delta_p(x))_i = \sum\limits_{j\in V}w_{ij}\phi_p(x_i-x_j),\;\; i\in\mathcal{I}_n
\end{equation}
Note that $\Phi_2$ is simply the identity function and for $p=2$, \eqref{eq:plapexpli} gives
\begin{equation}
(\Delta_2(x))_i = \sum\limits_{j\in V}w_{ij}(x_i-x_j) = (Lx)_i,\;\;i\in\mathcal{I}_n.
\end{equation}
We are now equipped with the necessary notation to define the $p$-Laplacian eigenproblem, which is the focus of the paper.
\begin{definition}[$p$-Laplacian eigenproblem]\label{def:plapprob}
Find $(\lambda,x)$ in $\mathbb{R}\times \mathbb{R}^n$ with $x\neq 0$ such that
\begin{equation}\label{eq:plapprob}
(\Delta_p(x))_i = \lambda \phi_p(x)_i,\;\; i\in\mathcal{I}_n.
\end{equation}
\end{definition}
The motivation for the definition of the eigenproblem as it  is presented in Definition~\ref{def:plapprob} can be found in Theorem~3.1 in \cite{Hein:2009:PLAPLACIAN}, which states that $x_*$ is a critical point of the functional
\begin{equation}
Q_p(x) = \frac{x^T\Delta_p(x)}{\norm{x}_p^p}
\end{equation}
if and only if $x_*$ is a solution to \eqref{eq:plapprob} with eigenvalue $\lambda_* = Q_p(x_*)$. For $p = 2$, \eqref{eq:plapprob} reduces to the linear eigenvalue problem for $L$.

\subsection{NEPv1 formulation}\label{sec:nepvform}
We now proceed to elaborate on the first of the two main ideas of the paper, which is to reformulate the $p$-Laplacian eigenproblem as a specific type of a problem that is known in the numerical linear algebra literature as the \textit{nonlinear eigenproblem with eigenvector nonlinearity} (NEPv). In the most general setting, a NEPv can be considered as a generalization of the invariant subspace problem for linear operators.
\begin{definition}[NEPv]\label{def:nepvdef}
Given $M:\mathbb{R}^{n\times r}\to\mathbb{R}^{n\times n}$, find $(\Lambda,X) \in \mathbb{R}^{r\times r}\times \mathbb{R}^{n\times r}$ where $\Lambda$ is a diagonal matrix and
\begin{equation}\label{eq:nepvdef}
\begin{split}
M(X)X &= X \Lambda,\\
X^TX &= I_p,
\end{split}
\end{equation}
where $M(XQ) = M(X)$ for any unitary $Q$.
\end{definition}
The columns of $X$ are called the eigenvectors of $M$ and the diagonal entries of $\Lambda$ are called the eigenvalues. Note that the first subequation in \eqref{eq:nepvdef} can be separated into $r$ different equations.
\begin{equation}
M(X)x_i = \lambda_i x_i,\;\; i=1,\ldots,r.
\end{equation}
where $\lambda_i$ is the $i$th diagonal element of $\Lambda$.

In the context of the $p$-Laplacian eigenproblem, a special case of the problem in Definition~\ref{def:nepvdef}, with $r = 1$ is relevant. We call it a \textit{nonlinear eigenproblem with \textbf{one} eigenvector nonlinearity}, or NEPv1.
\begin{definition}[NEPv1]\label{def:nepv1def}
Find $(\lambda,x)$ in $\mathbb{R}\times \mathbb{R}^{n}$ with $x\neq 0$ such that
\begin{equation}\label{eq:nepv1def}
M(x)x = \lambda x,
\end{equation}
where $M(\alpha x) = M(x)$ for any $\alpha\in \mathbb{R}$.
\end{definition}

We turn our attention back to the $p$-Laplacian eigenproblem and show how it can be reformulated in the form of a NEPv1. Consider $f:\mathbb{R}^n\to\mathbb{R}^m$ that maps a vector defined over the set of nodes to one that is defined over the edges, such that
\begin{equation}
f(x) = D^w\Phi_p(Bx).
\end{equation}
Hence, $(f(x))_k = w_{ij}\Phi_p(x_i-x_j)$ where $w_{ij}$ is the weight of edge $k$ with $k_1 = j$ and $k_2 = i$. In other words, for $x$ defined on the nodes, $(f(x))_k$ contains for edge $k$ its weight times $\Phi_p$ applied to the difference of $x_{k_2}$ and $x_{k_1}$. Left-multiplying $f(x)$ by $B^T$, we have
\begin{equation}\label{eq:plap_bdb}
(B^TD^w\Phi_p(Bx))_i = \sum\limits_{j\in V}w_{ij}\Phi_p(x_i-x_j),\;\; i=1,\ldots,n
\end{equation}
Comparing \eqref{eq:plapexpli} and \eqref{eq:plap_bdb}, we have
\begin{equation}
\Delta_p(x) = B^TD^w\Phi_p(Bx)
\end{equation}
and the $p$-Laplacian eigenproblem as specified by \eqref{eq:plapprob} can be written as
\begin{equation}\label{eq:plapprob_alt}
B^TD^w\Phi_p(Bx) = \lambda \Phi_p(x).
\end{equation}

We note that for any $y\in\mathbb{R}^l$ such that $y_i \neq 0$ for all $i\in \{1,\ldots,l\}$, we have
\begin{equation}\label{eq:phipexp}
\Phi_p(y) = \diag(|y|)^{p-1}\sign(y) = \diag(|y|)^{p-2}y,
\end{equation}
where $\diag(|y|)$ denotes the diagonal matrix containing aboslute values of the elements in $y$ as the diagonal entries.

Using \eqref{eq:phipexp} with $y = Bx$ in the left hand side of \eqref{eq:plapprob_alt} and $y = x$ in the right hand side of \eqref{eq:plapprob_alt} gives us
\begin{subequations}\label{eq:pform1}
\begin{eqnarray}
B^TD^w\diag(|Bx|)^{p-2}Bx &=& \lambda \diag(x)^{p-2}x,\label{eq:pform1a}\\
  (Bx)_i &\neq& 0,\;\;i=1,\ldots,m\label{eq:pform1b}\\
   x_i &\neq & 0,\;\;i=1,\ldots,n\label{eq:pform1c}
\end{eqnarray}
\end{subequations}
To ease notation, we define the operator $\mathcal{D}:\mathbb{R}^l\to\mathbb{R}^{l\times l}$ that takes as input a vector of arbitrary length and outputs the action of $\diag$ on the vector.
\begin{equation}
\mathcal{D}(y) = \diag(y).
\end{equation}

The first subequation \eqref{eq:pform1a} is in the form of a what we call a generalized NEPv1.
\begin{definition}[Generalized NEPv1]\label{def:gnepv1def}
Find $(\lambda,x) \in \mathbb{R}\times \mathbb{R}^{n}$ with $x\neq 0$ such that
\begin{equation}\label{eq:gnepv1def}
N(x)x = \lambda R(x)x,
\end{equation}
where there exists $g:\mathbb{R}\setminus\{0\}\to\mathbb{R}\setminus\{0\}$ such that
\begin{equation}
N(\alpha x) = g(\alpha)N(x),\;\;
R(\alpha x) = g(\alpha)R(x)
\end{equation}
for any $\alpha \neq 0$.
\end{definition}
For \eqref{eq:pform1}, we have
\begin{subequations}
\begin{eqnarray}
N(x) &=& B^TD^w\dd(|Bx|)^{p-2}B,\label{eq:nrdef1a}\\
R(x) &=& \dd(|x|)^{p-2}\label{eq:nrdef1b}.
\end{eqnarray}\label{eq:nrdef1}
\end{subequations}
and
\begin{equation}
g(\alpha) = |\alpha|^{p-2}.
\end{equation}
The matrix $R(x)$ is invertible for any $x$ satisfying \eqref{eq:pform1c} and hence, \eqref{eq:gnepv1def} can be rewritten by multiplying both sides by $R(x)^{-1}$, that is
\begin{equation}\label{eq:gtonepv1}
R(x)^{-1}N(x)x = \lambda x,
\end{equation}
which is a NEPv1 of the type defined by \eqref{eq:nepv1def} with $M(x) = R(x)^{-1}N(x)$. Although \eqref{eq:gtonepv1} and \eqref{eq:gnepv1def} are equivalent, the generalized NEPv1 formulation in \eqref{eq:gnepv1def} is more suitable because of the symmetry of $N$, which is an important property to have for successful convergence of the SCF iteration, as we explain in the Section~\ref{sec:scfplap}. However, before we proceed to Section~\ref{sec:scfplap}, we would like to introduce another approach to formulate the $p$-Laplacian eigenproblem as a generalized NEPv1.

We first notice that for any $y\in\mathbb{R},z \in \mathbb{R}$,
\begin{equation}\label{eq:sign}
\sign(y-z) = \begin{cases}
              0\;\;&\textrm{if}\;y = z,\\
              \sign(y)\;\textrm{or}\;-\sign(z),\;\;&\textrm{if}\;y = -z,\\
              \sign(y),\;\;&\textrm{if}\;|y| > |z|,\\
              -\sign(z),\;\;&\textrm{if}\;|y| < |z|.
         \end{cases}
\end{equation}

We define a function $P:\mathbb{R}^{n}\to\mathbb{R}^{m\times n}$ as follows.
\begin{equation}\label{eq:pdef}
P(x)_{k} = \begin{cases}
            B_{k},&\;\;\textrm{if}\;x_{k_1} = x_{k_2}\\
            e_{k_2}^T,&\;\;\textrm{if}\;x_{k_1} = -x_{k_2}\\
            e_{k_2}^T,&\;\;\textrm{if}\;\big|x_{k_1}\big| < \big|x_{k_2}\big|\\
            -e_{k_1}^T,&\;\;\textrm{if}\;\big|x_{k_1}\big| > \big|x_{k_2}\big|\\
            \end{cases}
\end{equation}
For any $x$ defined over the nodes, $P(x)$ is of the same size as $B$ and depends on the relationship between the absolute values of elements in $x$. The $k$th row of $P(x)$ depends only on $x_{k_1}$ and $x_{k_2}$. The way $P(x)$ has been defined in \eqref{eq:pdef} leads to a useful identity for all $x\in\mathbb{R}^{n}$, given by
\begin{equation}\label{eq:Pbid}
\sign(Bx) = P(x)\sign(x).
\end{equation}
This would also hold for some alternate definitions of $P(x)$. For example, we could define $P$ as follows.
\begin{equation}\label{eq:pdef2}
P(x)_{k} = \begin{cases}
            0,&\;\;\textrm{if}\;x_{k_1} = x_{k_2}\\
            e_{k_2}^T,&\;\;\textrm{if}\;x_{k_1} = -x_{k_2}\\
            e_{k_2}^T,&\;\;\textrm{if}\;\big|x_{k_1}\big| < \big|x_{k_2}\big|\\
            -e_{k_1}^T,&\;\;\textrm{if}\;\big|x_{k_1}\big| > \big|x_{k_2}\big|\\
            \end{cases}
\end{equation}
The identity \eqref{eq:Pbid} is useful as it allows us to formulate \eqref{eq:plapprob_alt} as
\begin{equation}\label{eq:preform2}
B^TD^w\dd\left(|Bx|\right)^{p-1}P(x)\sign(x) = \lambda \dd\left(|x|\right)^{p-1}\sign(x).
\end{equation}
This can be used to characterize a necessary condition for the existence of solutions to the $p$-Laplacian eigenproblem in terms of existence of \textit{signed eigenvectors} (eigenvectors whose elements can only be $-1$, $1$ or $0$) to a generalized eigenvalue problem.
\begin{theorem}[Necessary condition for solutions to \eqref{eq:plapprob}]
If the $p$-Laplacian eigenproblem defined by \eqref{eq:plapprob} has a solution $(\lambda_*,x_*)$, then the generalized eigenvalue problem defined as
\begin{equation}\label{eq:geneig_nec}
N'y = \lambda R'y
\end{equation}
with
\begin{equation}\label{eq:nprp}
N' = B^TD^w\dd\left(|Bx_*|\right)^{p-1}P(x_*),\;\; R' =\dd\left(|x_*|\right)^{p-1}
\end{equation}
must have a signed eigenvector $y_*$ with eigenvalue $\lambda_*$. Moreover, $y_* = \sign(x_*)$.
\end{theorem}
This motivates a theoretical investigation into the following question:\textit{Under what conditions on $N'$ and $R'$ does the generalized eigenvalue problem in \eqref{eq:geneig_nec} have signed eigenvectors?}. If such conditions can be identified and we find that \eqref{eq:nprp} does not satisfy those conditions, then \eqref{eq:plapprob} will not have a solution.

We note that \eqref{eq:preform2} is not yet in the form of a generalized NEPv1. For $x$ with nonzero elements, it can be rewritten as
\begin{subequations}
\begin{eqnarray}
B^TD^w\dd\left(|Bx|\right)^{p-1}P(x)\dd\left(|x|\right)^{-1}x &=& \lambda \dd\left(|x|\right)^{p-2}x,\label{eq:form2a}\\
x_i &\neq& 0,\;\;i=1,\ldots,n,\label{eq:form2b}
\end{eqnarray}\label{eq:form2}
\end{subequations}
which is a generalized NEPv1 of the type \eqref{eq:gnepv1def} with an additional constraint \eqref{eq:form2b} and $N$ and $R$ as defined by
\begin{subequations}
\begin{eqnarray}
N(x) &=& B^TD^w\dd\left(|Bx|\right)^{p-1}P(x)\dd\left(|x|\right)^{-1},\label{eq:nrdef2a}\\
R(x) &=& \dd\left(|x|\right)^{p-2}.\label{eq:nrdef2b}
\end{eqnarray}\label{eq:nrdef2}
\end{subequations}
and
\begin{equation}
g(\alpha) = |\alpha|^{p-1}.
\end{equation}
Until now, we have introduced two separate ways to formulate the $p$-Laplacian eigenproblem as a generalized NEPv1. We refer to the formulation of \eqref{eq:gnepv1def} with $N$ and $R$ as defined by \eqref{eq:nrdef1} and constrained by \eqref{eq:pform1a} and \eqref{eq:pform1b} as \textit{Form 1}. We call the formulation with $N$ and $R$ as defined by \eqref{eq:nrdef2} together with the constraint \eqref{eq:form2b} as \textit{Form 2}.

Both Form 1 and Form 2 contain constraints that were not present in the original formulation of the $p$-Laplacian eigenproblem given by \eqref{eq:plapprob_alt}. Form 1 has an additional constraint that is not present in Form 2. Hence, if $\mathcal{S}_0$, $\mathcal{S}_1$ and $\mathcal{S}_2$ denote the solution set to the original problem, Form 1 and Form 2 respectively, then we have $\mathcal{S}_1 \subseteq \mathcal{S}_2 \subseteq \mathcal{S}_0$.
\begin{table}[htbp!]
\begin{center}
\begin{tabular}{|c|c|c|c|c|c|}
\hline
& $N(x)$& $R(x)$& $x_i \neq 0$& $(Bx)_i \neq 0$& $N(x) = N(x)^T$\\
\hline
F1&  $B^TD^w\dd(|Bx|)^{p-2}B$& $\dd\left(|x|\right)^{p-2}$& Yes& Yes& Yes\\
\hline
F2& $B^TD^w\dd\left(|Bx|\right)^{p-1}P(x)\dd\left(|x|\right)^{-1}$& $\dd\left(|x|\right)^{p-2}$& Yes& No& No\\
\hline
\end{tabular}
\end{center}
\caption{Form 1 and Form 2 (F1 and F2), a comparison}
\label{tab:f1f2}
\end{table}

\section{The SCF iteration for the $p$-Laplacian eigenproblem}\label{sec:scfplap}
The SCF iteration is a fixed point iteration for solving a NEPv. The algorithm starts with an initial guess for the eigenvectors and in each step, solves a linear eigenvalue problem for a matrix that depends on the eigenvectors computed from the previous step. We illustrate the algorithm for a NEPv1 as defined in \eqref{eq:nepv1def} in Algorithm~\ref{alg:scf}.
\begin{algorithm}
\SetKw{Kw}{break}
\KwData{Intial guess: $v_0 \in \mathbb{R}^{n}$}
\KwResult{$(\lambda_*,v_*) \in \mathbb{R}\times \mathbb{R}^{n}$}
\For{$k = 0,1,\ldots$ until convergence}{
  Select appropriate $(\lambda_{k+1},v_{k+1})$ such that
  \[
  M(v_{k})v_{k+1} = \lambda_{k+1}v_{k+1}
  \]
}
$\lambda_* = \lambda_k$\\
$v_* = v_k$
\caption{The self-consistent field (SCF) iteration for a NEPv1}
\label{alg:scf}
\end{algorithm}

Note that in each step $k$ inside the for loop in Algorithm~\ref{alg:scf}, we have to select only one among the many eigenpairs of $M(v_{k})$. The selection strategy ususally depends on the context provided by the application. For example, for the solution of the Gross-Pitaevskii equation in quantum physics \cite{IMPLICIT:2021}, we would select the eigenpair corresponding to the smallest eigenvalue because we are interested in computing the ground state of bosons. For this paper, we will focus on applications that require computing the eigenpair corresponding to the second smallest eigenvalue, that is, a partition into two clusters. In practical applications, this approach is particularly useful if the dataset is expected to have a recursive structure \cite{YU:2020:PECOS}. Then we can use a hierarchical clustering technique where in each level or generation, the new clusters are obtained by splitting the clusters from the previous generation into two.

For applications in spectral clustering, we are exclusively interested in computing the real eigenpairs of the $p$-Laplacian operator. This means that in each step of the SCF iteration, we want $(\lambda_{k+1},v_{k+1})$ to be real, which is guaranteed if $M(v_k)$ in Algorithm~\ref{alg:scf} is symmetric. Moreover, there are several convergence results for the SCF iteration \cite{DENSITY:2021,Cances:2000:SCF,Bai:2020:OPTIMALSCF} that characterize its behaviour and conditions for convergence in the case of symmetric (Hermitian for the complex case) $M$. Similar results are not available for the non-symmetric case. Hence, in order to apply the SCF iteration to solve the generalzed NEPv1 formulations derived in the previous section, we cast them in a symmetric NEPv1 form.

For Form 1, with the assumption that the constraint \eqref{eq:pform1c} is satisfied by $x$, we do this by rewriting \eqref{eq:gnepv1def} as
\begin{equation}\label{eq:gnepv1scf}
R(x)^{\sm 1/2}N(x)R(x)^{\sm 1/2}y(x) = \lambda y(x),
\end{equation}
where
\begin{equation}
y(x) = R(x)^{1/2}x.
\end{equation}
We note that $R(x)$ is a diagonal matrix and $N(x) = B^TD^w\dd(|Bx|)^{p-2}B$ is symmetric, which means $R(x)^{\sm 1/2}N(x)R(x)^{\sm 1/2}$ is symmetric. We use \eqref{eq:gnepv1scf} to adapt the SCF iteration to solve Form 1, and present this adaptation of Algorithm~\ref{alg:scf} as Algorithm~\ref{alg:scfgnepv1}.
\begin{algorithm}
\SetKw{Kw}{break}
\KwData{Intial guess: $v_0 \in \mathbb{R}^{n}$}
\KwResult{$(\lambda_*,v_*) \in \mathbb{R}\times \mathbb{R}^{n}$}
\For{$k = 0,1,\ldots$ until convergence}{
  Select appropriate $(\lambda_{k+1},y_{k+1})$ such that
  \[
  R(v_k)^{\sm 1/2}N(v_{k})R(v_k)^{\sm 1/2}y_{k+1} = \lambda_{k+1}y_{k+1}
  \]
  Transform $y_{k+1}$ to $v_{k+1}$ by
  \[
   v_{k+1} = R(v_k)^{\sm 1/2}y_{k+1}
  \]
}
$\lambda_* = \lambda_k$\\
$v_* = v_k$
\caption{The SCF iteration for Form 1 with $N$, $R$ given by \eqref{eq:nrdef1} and the assumption that the constraints $\eqref{eq:pform1b}$ and $\eqref{eq:pform1c}$ is satisfied by the iterates $v_k$.}
\label{alg:scfgnepv1}
\end{algorithm}

In the case of Form 2, $N$ is not symmetric and hence, it is not obvious how to obtain a symmetric reformulation like in \eqref{eq:gnepv1scf}. Moreover, computing $N(x)$ involves computing $P(x)$, which requires $\mathcal{O}(m)$ operations in each step and can be prohibitively expensive to compute in each step of the SCF iteration. Also, by definition, $P$ is not a continuous function of $x$. Hence, although Form 2 contains fewer constraints than Form 1 and provides a tool for further theoretical investigation, it is not suited for applying the SCF iteration, the successful convergence of which depends on the symmetry and smoothness of $N$ \cite{DENSITY:2021}.
\subsection{Regularization by smooth approximation of the absolute value function}
We have adapted the SCF iteration for solving the $p$-Laplacian eigenproblem using Form 1 in Algorithm~\ref{alg:scfgnepv1}. However, we have assumed in the algorithm that the constraints of Form 1, namely \eqref{eq:pform1b} and \eqref{eq:pform1c} are always satisfied by the iterates $v_k$. This need not be the case, particularly when $v_k$ is close to the target eigenvector $v_*$, which in our case is the second eigenvector, meaning the one corresponding to the second smallest eigenvalue.

The second eigenvector $v_*$ of the $p$-Laplacian is known to contain elements which are close at indices corresponding to the nodes which may belong in the same cluster. This becomes even more noticeable the closer we get to the limit $p \to 1^+$ and leads to a natural technique to recover the two clusters from the second eigenvector by applying a suitable threshold that separates the entries into their respective clusters \cite{Hein:2009:PLAPLACIAN}. However, this also implies that when $v_k$ is close to convergence towards $v_*$, the entries of $v_k$ will display similar behaviour meaning that for some edges $j$, $(Bv_{k})_j$ will be close to zero. These leads to numerical issues with ill-conditioned matrices because some entries in $N(v_k) = B^TD^w\dd(|Bx|)^{p-2}B$ become quite large as $v_k\to v_*$.

To alleviate such issues, we modify Form 1 by approximating the absolute value function with another one that behaves similar to it for nonzero inputs, but is differentiable and positive at zero. There are several papers in the optimization and machine learning literature for accurate and smooth approximations of the absolute value function \cite{VORONIN:2015,SCHMIDT:2007,BRIDLE:SOFTMAX:1989}. We expand $|x|$ as
\begin{equation}\label{eq:modx}
|x| = \max(x,0)+\max(-x,0)
\end{equation}
and then apply the approximation used in \cite{CHEN:1996} for the $\max$ function, which can be obtained by integrating the the sigmoid function given by $\operatorname{sigmoid}(x) = 1/(1+e^{-ax})$. This approach gives the approximation
\begin{equation}\label{eq:max}
\max(x,0) \approx x+\frac{\log(1+e^{-ax})}{a},
\end{equation}
where $a$ is a parameter that can be used to control the smoothness and accuracy of the approximation. Combining \eqref{eq:max} and \eqref{eq:modx}, we have an approxmation for the absolute value function given by
\begin{equation}\label{eq:softabs}
|x| \approx \frac{\log(1+e^{-ax})+\log(1+e^{-ax})}{a} := \operatorname{sf}_{a}(x),
\end{equation}
which we refer to as the \textit{softabs} fucntion, similar in name to the \textit{softmax} \cite{BRIDLE:SOFTMAX:1989} and \textit{softplus} \cite{SOFTPLUS:2000} functions, which are ubiquitous in machine learning. We generalize the softabs function to act elementwise on vector-valued inputs $y\in\mathbb{R}^l$ as
\begin{equation}
\operatorname{sf}_a(y)_i = \frac{\log(1+e^{-ay_i})+\log(1+e^{-ay_i})}{a},\;\;i=1,\ldots,n.
\end{equation}
\begin{figure}[h]
\overfullrule=0pt
  \begin{subfigure}[b]{0.5\linewidth}
  \centering
  \includegraphics[scale=0.7]{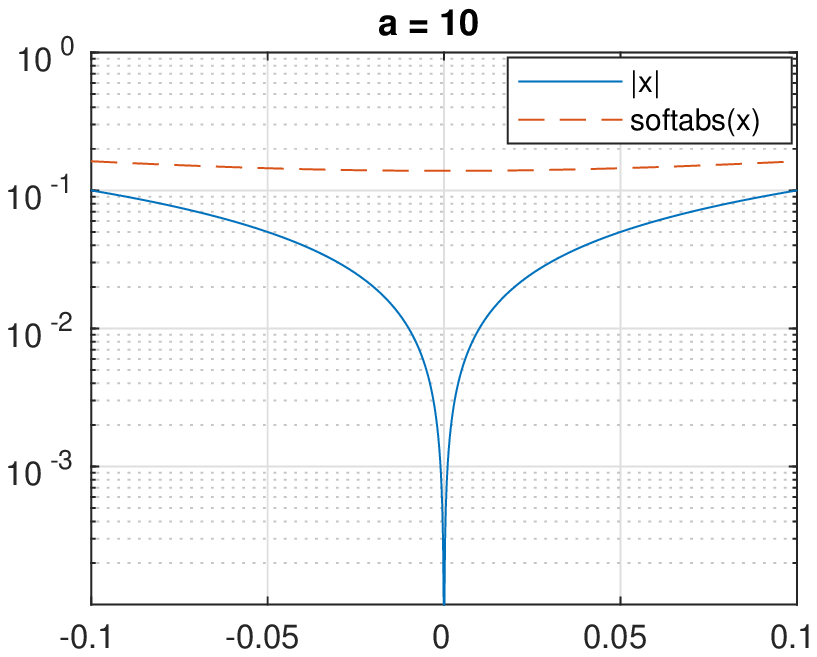}
  \end{subfigure}
  \begin{subfigure}[b]{0.5\linewidth}
  \centering
  \includegraphics[scale=0.7]{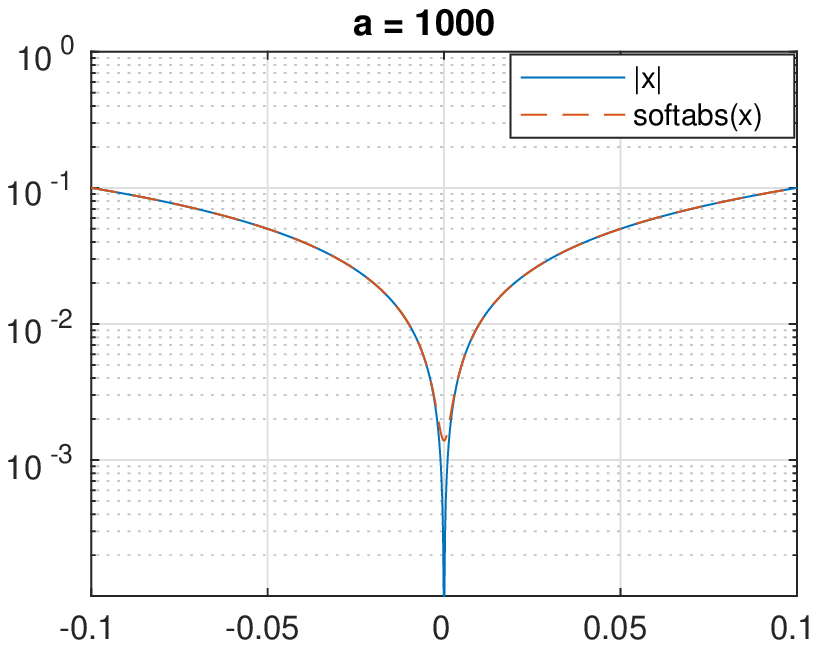}
  \end{subfigure}
  \caption{A comparison of $\operatorname{sf}_a(x)$ and $|x|$ for two different values of $a$, plotted with a log scale for the vertical axis}
  \label{fig:softabs}
\end{figure}
Note that $sf_a(x)\to |x|$ pointwise as $a\to\infty$. In particular, for any $x\in\mathbb{R}$, we have
\begin{equation}
|\operatorname{sf}_a(x)-|x|| \leq sf_a(0) = \mathcal{O}(1/a).
\end{equation}
We now consider the \textit{regularized Form 1} with $N$ and $R$ modified to $\tilde{N}_a$ and $\tilde{R}_a$ obtained by applying the softabs approximation \eqref{eq:softabs}.
\begin{subequations}
\begin{eqnarray}
\tilde{N}_a(x) &=& B^TD^w\dd\left(\operatorname{sf}_a(Bx)\right)^{p-2}B,\label{eq:nrtdefa}\\
\tilde{R}_a(x) &=& \dd\left(\operatorname{sf}_a(x)\right)^{p-2},\label{eq:nrtdefb}
\end{eqnarray}\label{eq:nrtdef}
\end{subequations}
Let $(\lambda_*,x_*)$ be an eigenpair of the modified problem, such that
\begin{equation}
\tilde{N}_a(x_*)x_* = \lambda_* \tilde{R}_a(x_*)x_*
\end{equation}
holds. Now, for Form 1 with $N$ and $R$ as defined by \eqref{eq:nrdef1}, $\Delta N(x_*) = N(x_*)-\tilde{N}_a(x_*)$ and $\Delta R(x_*) = R(x_*)-\tilde{R}_a(x_*)$ , we have
\begin{subequations}
\begin{eqnarray}
\norm{\Delta N(x_*)} &=& \norm{B^TD^w\left(\dd(|Bx|)^{p-2}-\dd(\operatorname{sf}_a(Bx))^{p-2}\right)B} \leq \mathcal{O}(1/a),\label{eq:deltanra}\\
\norm{\Delta R(x_*)} &=& \norm{\dd(|x|)^{p-2}-\dd(\operatorname{sf}_a(x))^{p-2}} \leq \mathcal{O}(1/a).\label{eq:deltanrb}
\end{eqnarray}\label{eq:deltanr}
\end{subequations}
Using \eqref{eq:deltanr}, we have
\begin{equation}
\begin{split}
\norm{\left(N(x_*)-\lambda_* R(x_*)\right)x_*} &\leq \norm{(\tilde{N}_a(x_*)-\lambda_* \tilde{R}_a(x_*))x_*+(\Delta N(x_*)-\lambda\Delta R(x_*))x_*}\\
 &\leq (\norm{\Delta N(x_*)}+\lambda_*\norm{\Delta R(x_*)})\norm{x_*}\\
 &\leq \mathcal{O}(1/a),
\end{split}
\end{equation}
which implies that any eigenpair of the regularized Form 1 is an approximate eigenpair of Form 1 with $N$ and $R$ defined by \eqref{eq:nrdef1} with residual norm bounded by $\mathcal{O}(1/a)$.

In applictions in data science and machine learning, one is often not interested in computing solutions to an accuracy of machine precision. In our context, using $a$ of the order of $10^5$ and computing the second eigenvector upto around five digits of accuracy is often enough. This is becuase the elements of the eigenvector corresponding to different clusters reveal themselves by a clear separation in their values and the clusters can be obtained by applying a simple threshold, as discussed before. Computing the eigenvector more accurately beyond this point does not necessarily provide better clustering.
\subsection{Intial guess and a $p$-iterated scheme}
Our experiments revealed that the selection of an appropriate initial guess for the eigenvector is important for the successful convergence of the SCF iteration applied to the regularized version of Form~1. The simplest candidate for an initial guess is the second eigenvector of the $2$-Laplacian, which is relatively easier to obtain computationally because it involves solving a linear eigenvalue problem for a matrix that is sparse and positive semi-definite.

However, if we to compute eigenpairs of the $p_{\textrm{target}}$-Laplacian where $p_{\textrm{target}}$ is close to one, then, the second eigenvector of the $2$-Laplacian can still be far from the target eigenvector. Hence, we use a scheme where we successively solve the regularized Form~1 for different values of $p$ using the SCF iteration, starting with $p = 2$ and in each step reduce $p$ by a suitable amount $\Delta p$ until we reach $p = p_{\textrm{target}}$. This $p$-iterated adaptation of the SCF iteration to solve the regularized Form 1 is presented in Algorithm~\ref{alg:scfsoft}.
\begin{algorithm}[htbp!]
\SetKw{Kw}{break}
\KwData{Target $p$: $p_{\textrm{target}}$, Change in $p$: $\Delta p$, Softabs parameter: $a$}
\KwResult{$(\lambda_*,v_*) \in \mathbb{R}\times \mathbb{R}^{n}$}
Solve for the second smallest eigenvector of the $2$-Laplacian to get $v_0$\\
$p_i = 2$\\
\While{$p_i \geq p_{\textrm{target}}$}{
 $p_i = p_i-\Delta p$\\
 Redefine $\tilde{N}_a$ and $\tilde{R}_a$ corresponding to the $p_i$-Laplacian\\
\For{$k = 0,1,\ldots$ until convergence}{
  Select appropriate $(\lambda_{k+1},v_{k+1})$ such that
  \[
  \tilde{R}_a(v_k)^{\sm 1/2}\tilde{N}_a(v_{k})\tilde{R}_a(v_k)^{\sm 1/2}y_{k+1} = \lambda_{k+1}y_{k+1}
  \]

  Transform $y_{k+1}$ to $v_{k+1}$ by
  \[
   v_{k+1} = \tilde{R}_a(v_k)^{\sm 1/2}y_{k+1}.
  \]
}
$v_0 = v_k$
}
$\lambda_* = \lambda_k$\\
$v_* = v_k$
\caption{The SCF iteration combined with the $p$-iterated scheme for the regularized Form 1 with $\tilde{N}_a$, $\tilde{R}_a$ given by \eqref{eq:nrtdef}}
\label{alg:scfsoft}
\end{algorithm}
\\
\section{Numerical examples}\label{sec:numex}
\subsection{Stochastic block model}
We consider a randomly generated graph containing $n$ nodes partitioned into two disjoint sets with $n_1$ and $n_2$ nodes respectively. The edges are generated such that there is a higher probability of having an edge between two randomly picked nodes if they belong to the same set compared to if they come from two different sets. Technically, this is called a \textit{stochastic block model} (SBM)\cite{HOLLAND:1983} with two communities, which we call $C_1$ and $C_2$. The adjacency matrix of the graph has a $2\times 2$ block structure given by
\begin{equation}
A = \begin{pmatrix}A_{11}& A_{12}\\A_{12}^T& A_{22}\end{pmatrix},
\end{equation}
where $A_{11}\in\mathbb{R}^{n_1\times n_1}$ and $A_{22} \in\mathbb{R}^{n_2\times n_2}$ are symmetric and can be seen as the adjacency matrices of the subgraphs obtained from $C_1$ and $C_2$ taken separately, respectively. The edges between the nodes in the two subgraphs are given in $A_{12}\in\mathbb{R}^{n_1\times n_2}$. The sparsity of these submatrices is encoded by a symmetric $2\times 2$ matrix
\begin{equation}
Q = \begin{pmatrix}q_{\textrm{in},1}& q_{\textrm{out}}\\ q_{\textrm{out}}& q_{\textrm{in},2}\end{pmatrix}
\end{equation}
where the $q_{\textrm{in},i}$ denotes the probability of two nodes being connected by an edge when both belong to $C_i$, for $i=1,2$. The variable $q_{\textrm{out}}$ denotes the probability of an edge connecting a node in $C_1$ to a node in $C_2$.Hence, if $q_{\textrm{in},1}$ and $q_{\textrm{in},2}$ are relatively much larger than $q_{\textrm{out}}$, then the clusters $C_1$ and $C_2$ are well-separated.

For our experiments with the SBM, we let $n_1 = n_2 = n_c$, that is we generate a graph with $n = 2n_{c}$ nodes in total. The weights of the edges are all ones and $q_{\textrm{in},1} = q_{\textrm{in},2} = q_{\textrm{in}}$.
\begin{figure}[h]
\overfullrule=0pt
  \begin{subfigure}[b]{0.5\linewidth}
  \centering
  \includegraphics[scale=0.8]{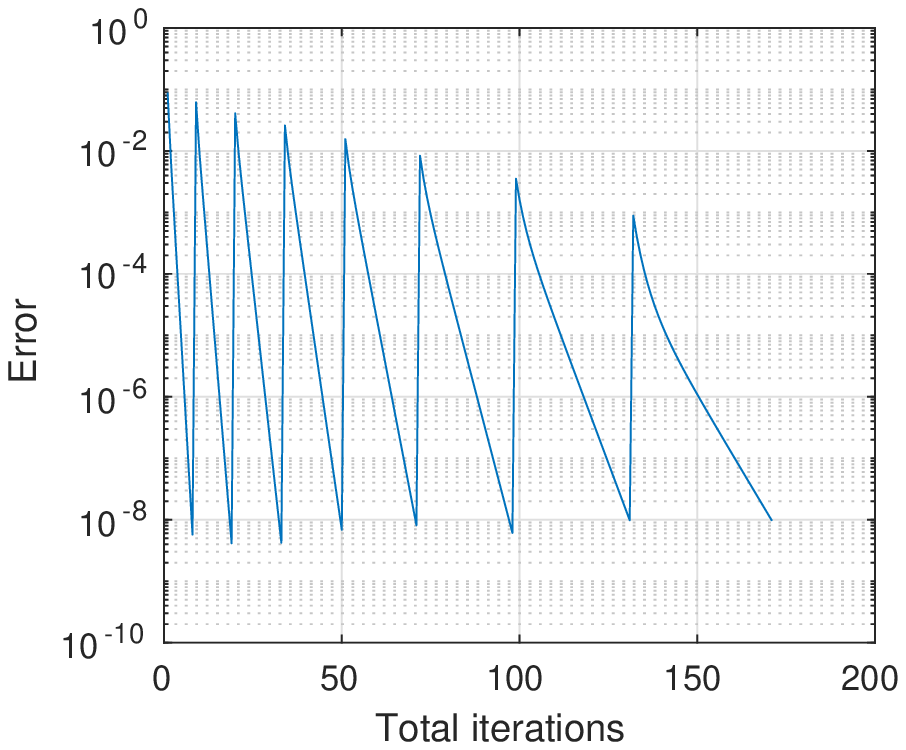}
  \caption{Error of the SCF iteration measured as $\norm{v_{k+1}-v_k}$}
  \label{fig:sbma}
  \end{subfigure}
  \begin{subfigure}[b]{0.5\linewidth}
  \centering
  \includegraphics[scale=0.8]{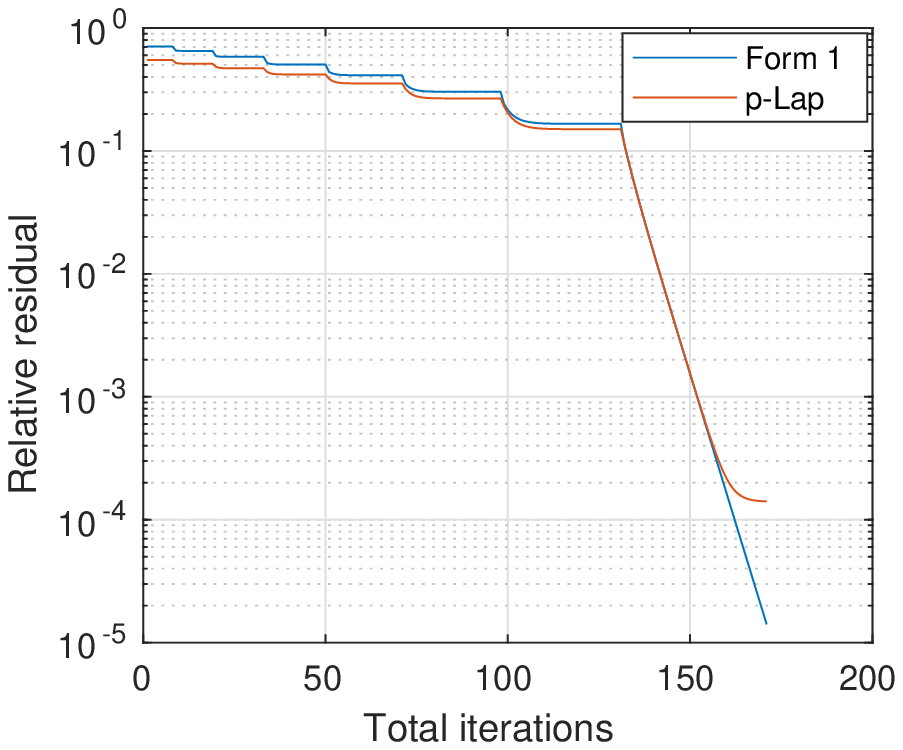}
  \caption{Relative residual of the regularized Form 1 and the $p$-Laplacian eigenproblem}
  \label{fig:sbmb}
  \end{subfigure}
  \caption{SCF iteration applied to the SBM model with $n_c = 100$, $q_{\textrm{in}} = 0.8$, $q_{\textrm{out}} = 0.3$, $a = 10^{10}$, $p_{\textrm{target}} = 1.2$}
  \label{fig:sbm_err_relres}
\end{figure}
Figure~\ref{fig:sbm_err_relres} shows convergence plots obtained by applying Algorithm~\ref{alg:scfsoft} with $p_{\textrm{target}} = 1.2$, $a = 10^{10}$, $\Delta p = 0.1$, to the problem generated by our SBM model with a total of two hundred nodes with $q_\textrm{in} = 0.8$, $q_\textrm{out} = 0.3$. The vertical axis in Figure~\ref{fig:sbma} shows the error of the SCF iteration, which at step $k$ measured by $\norm{v_{k+1}-v_k}$, the norm of the difference between the current and the future normalized eigenvector iterate. Figure~\ref{fig:sbmb} shows a plot of the relative residual of the $p$-Laplacian eigenproblem \eqref{eq:plapprob} at step $k$, measured using \eqref{eq:plapprob_alt} as $\norm{B^TD^{w}\Phi_{p_
{\textrm{target}}}(Bv_{k})-\lambda_k\Phi_{p_{\textrm{target}}}(v_k)}\big/{\norm{B^TD^{w}\Phi_{p_{\textrm{target}}}(Bv_{k})}}$, along with the relative residual of regularized Form 1, that is $\norm{\tilde{N}_a(v_k)v_k-\lambda_k \tilde{R}_a(v_k)}\big/\norm{\tilde{N}_a(v_k)v_k}$.

In Figure~\ref{fig:sbma}, we see successive SCF iterations converging one after another for each intermediate value of $p_i$ in Algorithm~\ref{alg:scfsoft}, until we reach $p_{\textrm{target}}$. The spikes are seen because the value of $p_i$ and the definition of $\tilde{N}_a$ and $\tilde{R}_a$ are updated upon convergence and the eigenvector that the iteration converged to becomes an initial guess for the new updated problem. Similar spikes are not observed in Figure~\ref{fig:sbmb} because the residual is always measured for the problem corresponding to $p = p_\textrm{target}$.

\subsubsection{Gaps between eigenvalues and convergence}
As we see in Figure~\ref{fig:sbm_err_relres}, the convergence of the SCF iteration is linear for each $p_i$ and the number of iterations required for convergence increases we reduce $p_i$. It is well known from the literature on its convergence analysis \cite{Cances:2000:SCF,Bai:2020:OPTIMALSCF,DENSITY:2021,Yang:2009:SCF} that the SCF iteration converges linearly, when it does converge. These results also bound the convergence factor $\textrm{c}_{\textrm{scf}}$ from above in terms of the gaps between the eigenvalues of the NEPv or NEPv1.

Athough tight upper bounds have been derived in some papers \cite{Bai:2020:OPTIMALSCF, DENSITY:2021}, the slowing down of the convergence with increasing $p_i$ in Figure~\ref{fig:sbma} can be explained using a simple conservative upper bound from \cite{Yang:2009:SCF} which bounds the convergence factor by an expression that is inversely proportional to the smallest gap $\delta$ between the so-called \textit{occupied} (the ones we are required want to compute) and \textit{unoccupied} (the ones that we donot require to compute) sets of eigenvalues, given by
\begin{equation}\label{eq:bound}
\textrm{c}_{\textrm{scf}} \leq \frac{\norm{K}}{\delta},
\end{equation}
where $K$ is a matrix depending on $\tilde{N}_a$ and $\tilde{R}_a$. In our case, the occupied set is a singleton set containing the second smallest eigenvalue, that is $\{\lambda_2\}$ and the unoccupied set contains all the other eigenvalues, that is $\{\lambda_1,\lambda_3,\ldots,\lambda_n\}$. Hence, we have
\begin{equation}\label{eq:delt}
\delta = \min\left(|\lambda_2-\lambda_1|,|\lambda_2-\lambda_3|\right).
\end{equation}
We measure $\delta$ as given by \eqref{eq:delt} for each $p_i$ upon convergence and present the results in Table~\ref{tab:gap}. From $p_i = 1.8$ onwards, $\delta$ goes down steadily, which increases the upper bound in \eqref{eq:bound} and hence, we expect the rate of convergence to steadily become slower as $p_i\to 1^+$.
\begin{table}
\centering
\begin{tabular}{|c|c|c|c|c|c|c|c|c|}
\hline
$p_i$& 1.9& 1.8& 1.7& 1.6& 1.5& 1.4& 1.3& 1.2\\
\hline
$\delta$&  0.0063& 0.0080& 0.0067& 0.0044& 0.0021& 0.0007& 0.0001& 0.000002\\
\hline
\end{tabular}
\caption{The gap $\delta$ for the regularized Form 1 specified by \eqref{eq:nrtdef} corresponding to different $p_i$ \label{tab:gap}}
\end{table}
\subsection{Two moons}
\begin{figure}[htbp!]
\overfullrule=0pt
  \begin{subfigure}[b]{0.5\linewidth}
  \centering
  \includegraphics[scale=0.8]{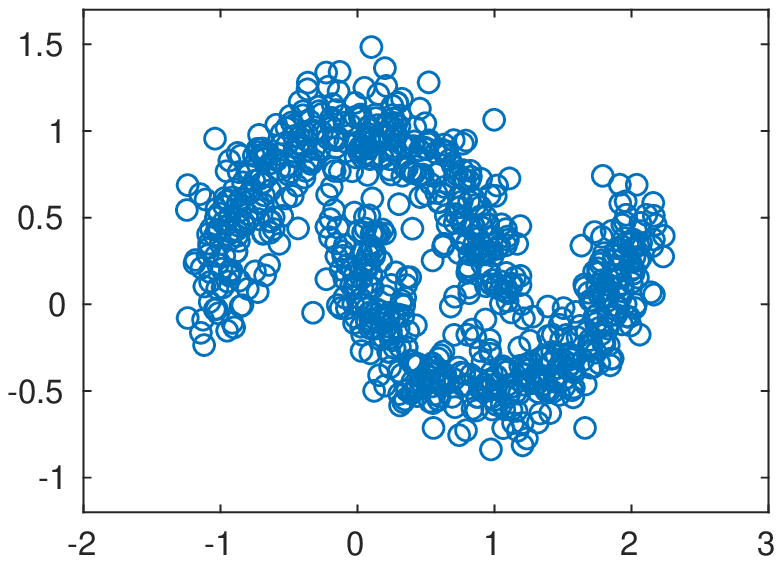}
  \caption{Projection in two dimensions \label{fig:2d}}
  \end{subfigure}
  \begin{subfigure}[b]{0.5\linewidth}
  \centering
  \includegraphics[scale=0.8]{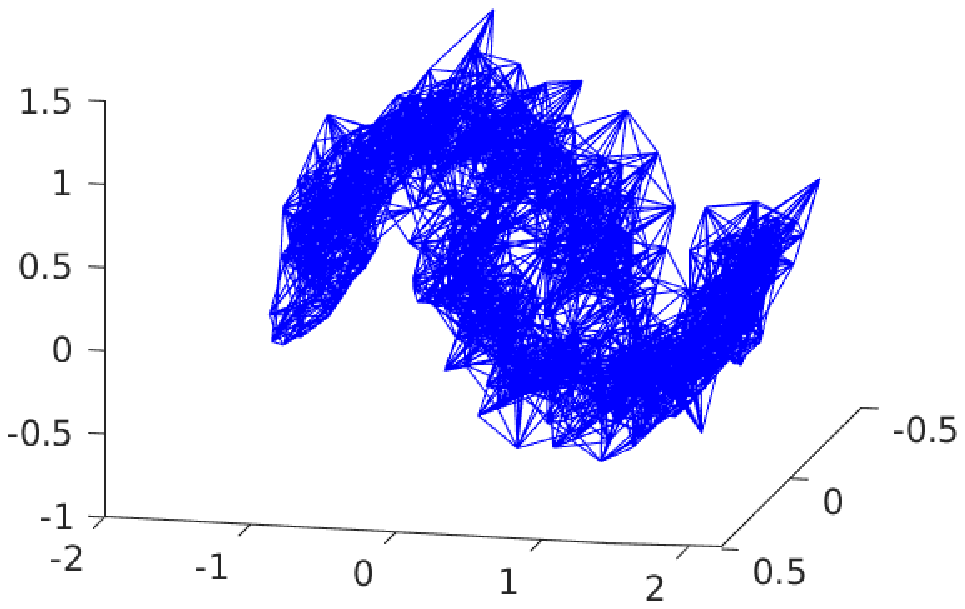}
  \caption{Three-dimensional projection with edge structure \label{fid:3d}}
  \end{subfigure}
  \caption{The two moons dataset \label{fig:twomoonsprojs}}
\end{figure}

\begin{figure}[htbp!]
\overfullrule=0pt
  \begin{subfigure}[b]{0.5\linewidth}
  \centering
  \includegraphics[scale=0.8]{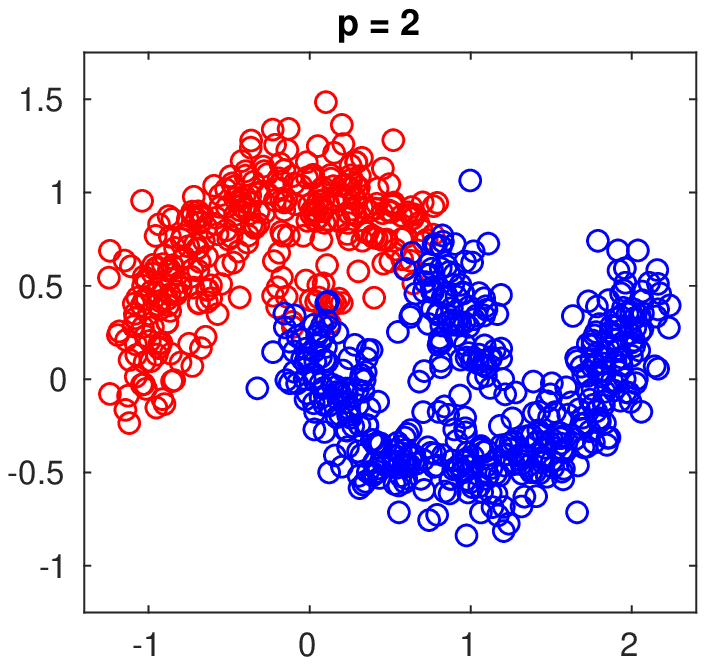}
  \end{subfigure}
  \begin{subfigure}[b]{0.5\linewidth}
  \centering
  \includegraphics[scale=0.8]{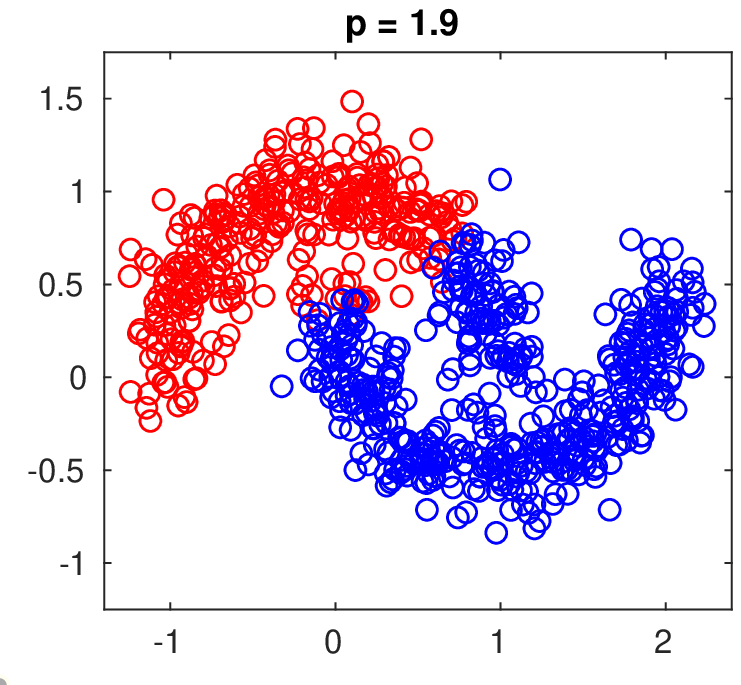}
  \end{subfigure}

  \begin{subfigure}[b]{0.5\linewidth}
  \centering
  \includegraphics[scale=0.8]{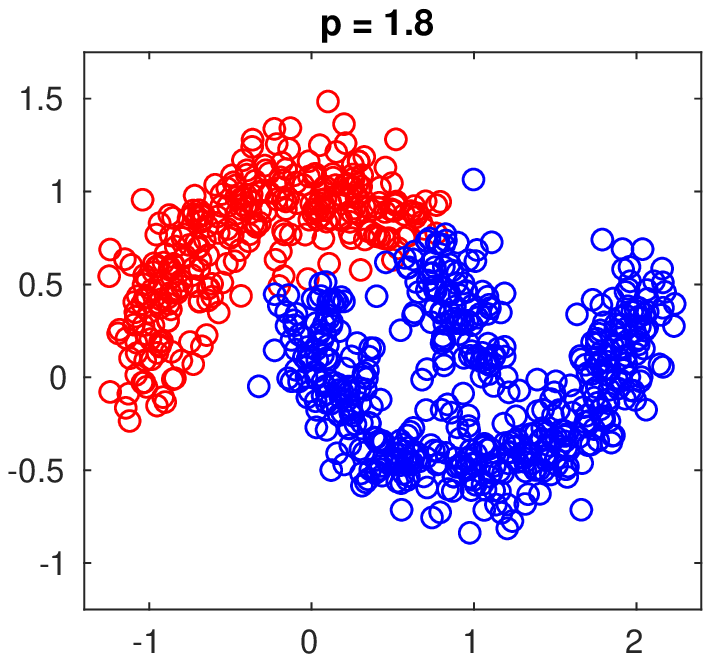}
  \end{subfigure}
  \begin{subfigure}[b]{0.5\linewidth}
  \centering
  \includegraphics[scale=0.8]{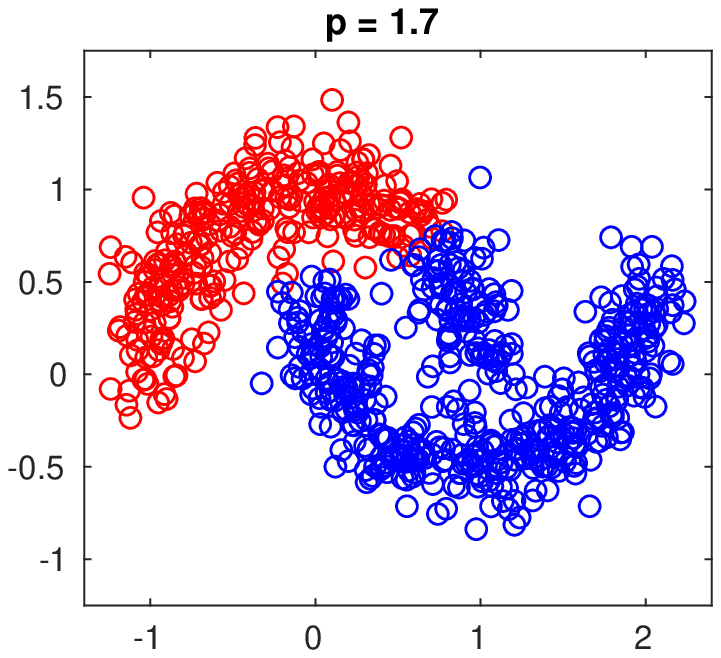}
  \end{subfigure}

  \begin{subfigure}[b]{0.5\linewidth}
  \centering
  \includegraphics[scale=0.8]{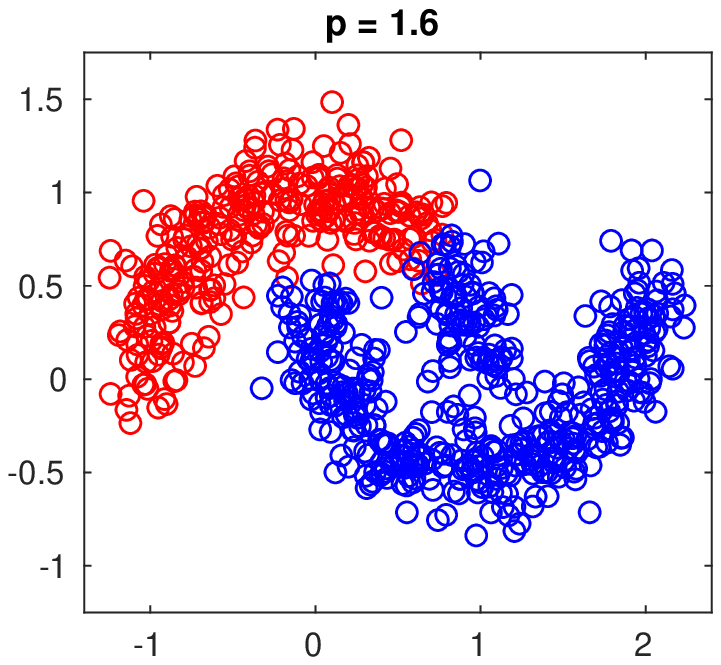}
  \end{subfigure}
  \begin{subfigure}[b]{0.5\linewidth}
  \centering
  \includegraphics[scale=0.8]{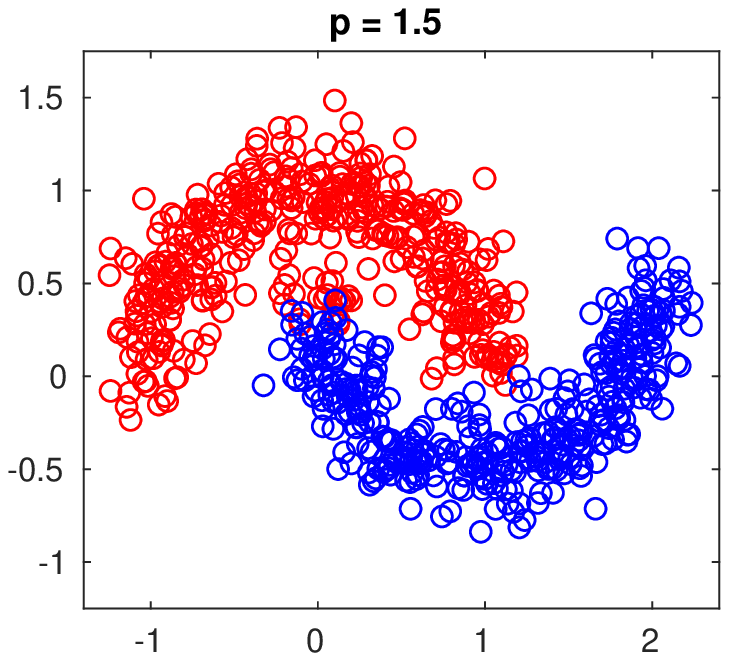}
  \end{subfigure}

  \begin{subfigure}[b]{0.5\linewidth}
  \centering
  \includegraphics[scale=0.8]{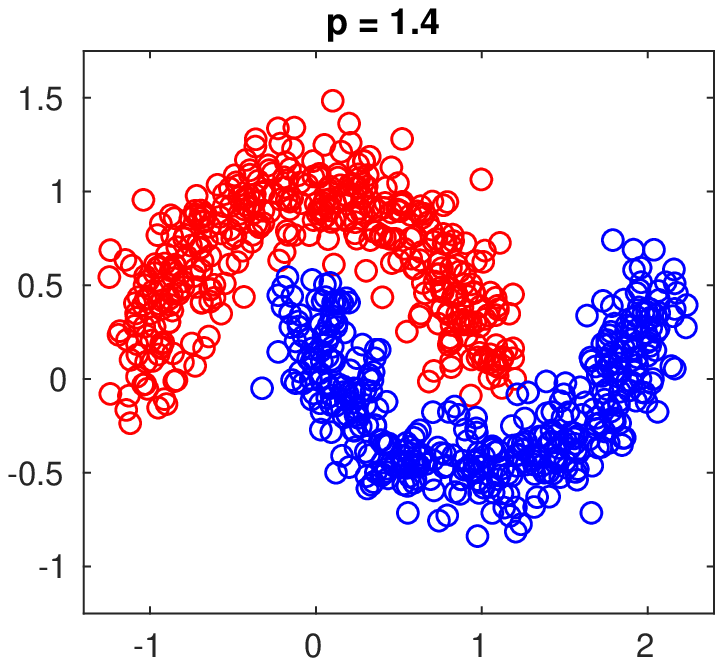}
  \end{subfigure}
  \begin{subfigure}[b]{0.5\linewidth}
  \centering
  \includegraphics[scale=0.8]{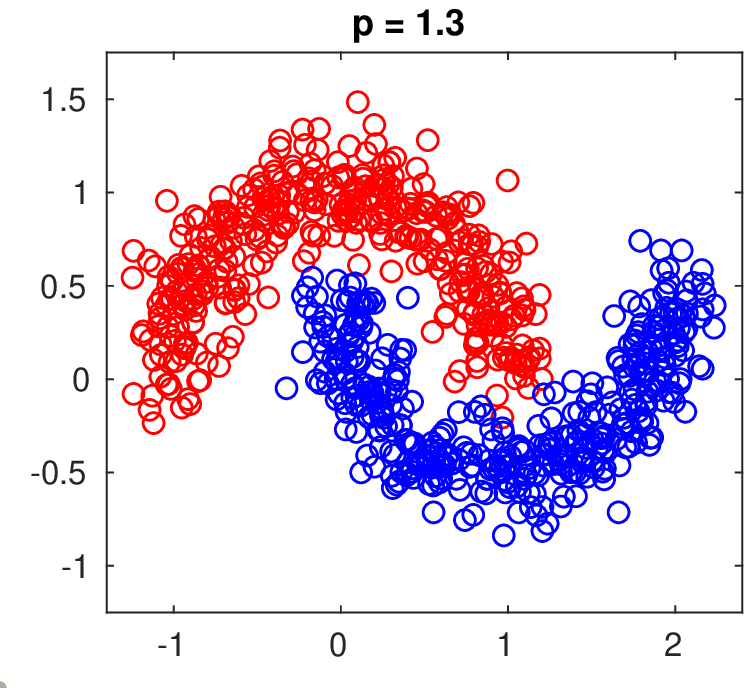}
  \end{subfigure}
  \caption{Clusters obtained in the two moons dataset by thresholding the second eigenvector of the $p$-Laplacian, for different $p$ \label{fig:twomoons}}
\end{figure}

The two moons dataset is created by generating points in two dimensions along two semicircles that are close to each other, and then adding a high dimensional noise to those points, as seen in Figure~\ref{fig:twomoonsprojs}. An ideal clustering of this dataset will have two clusters with each cluster corresponding to one of the semicircles. The addition of noise to the dataset corrupts the ideal structure of the dataset, and hence it serves as a good test for clustering algorithms. We consider a total of $n = 2n_c = 800$ nodes and set up the experiment in a similar manner as in \cite{Hein:2009:PLAPLACIAN}, by considering a Gaussian similarity function between the points where
\begin{equation}\label{eq:tmsim}
s(i,j) = \max(e^{\shortminus 2\norm{x_i-x_j}^2/\sigma_i^2},e^{\shortminus 2\norm{x_i-x_j}^2/\sigma_j^2}),
\end{equation}
where $\sigma_i$ is the euclidean distance between node $i$ and its nearest neighbour. We use \eqref{eq:tmsim} to build a symmetric $k$-nearest neighbours graph with $k = 10$. The two semicircles are chosen to be centred at $(0,0)$ and $(1,0.5)$ with radius $r = 1$. We add a $d$-dimensional Gaussian noise $\mathcal{N}(0,\sigma^2\mathbbm{1}_d)$ with $d = 10$ and $\sigma^2 = 0.02$.
The clusters are obtained from the second eigenvector by sorting the entries of the eigenvector, sweeping through all the different possible separation of entries by applying different thresholds and picking the threshold that minimizes RCC, as given by \eqref{eq:rcc}. This is done for every $p_i$ in Algorithm~\ref{alg:scfsoft}, and the resulting clusters are shown in Figure~\ref{fig:twomoons}.
As we see from Figure~\ref{fig:twomoons}, the clusters obtained for higher values of $p>1.5$ appear to be of a lower quality compared those obtained for $p\leq 1.5$. This is expected theoretically since the optimal values of RCC are obtained as $p\to 1$. This theoretical expectation is further confirmed by the data in Table~\ref{tab:rcc}, where we see that both the ratio cut and the ratio Cheeger cut decrease as we decrease $p$.
\begin{table}
  \centering
  \begin{tabular}{|c|c|c|c|c|c|}
  \hline
  $p_i$& 2.0& 1.8& 1.6& 1.4& 1.2\\
  \hline
  RCC& 0.02630& 0.02488& 0.02190& 0.01362& 0.01309\\
  \hline
  RCut& 0.04063&  0.03958& 0.03522& 0.02711& 0.02611\\
  \hline
  \end{tabular}
  \caption{RCC and RCut for different $p_i$\label{tab:rcc}}
\end{table}

\section{Conclusions and possible improvements}\label{sec:conc}
We have adapted the SCF iteration and shown that it can be used successfully to solve the $p$-Laplacian eigenproblem, as illustrated by the numerical experiments in the previous section. This provide a new bridge between the fields of quantum chemistry and data science. The formulation as a NEPv1 and the application of the SCF iteration to this problem was previously unexplored, and our results point towards several promising directions.

Our experiments haven been able to provide good results for upto $p=1.2$, but we are limited by the ill-conditioning of the matrices in Form 1 when $p$ is very close to one. This problem arises because entries in the eigenvector corresponding to nodes belonging to the same cluster become almost equal as $p\to 1$. A possible solution is to use intermediate agglomeration step whenever such a situation arises during the SCF iteration, that is to combine such nodes into a bigger supernode and modify the adjacency and incidence matrices accordingly. After the agglomeration step, we can continue the SCF iteration for the modified $p$-Laplacian of the smaller graph. The assumption contained in this approach is that nodes combined into a bigger supernode belong to the same cluster.

Although we focus only on computing the second eigenpair of the $p$-Laplacian, our approach can possibly be extended to compute several eigenpairs as well. In that case, the problem will not be a NEPv1, but a NEPv of the type defined in \eqref{eq:nepvdef}. One challenge in this direction is that the eigenvectors of the $p$-Laplacian are not orthogonal in general, but \eqref{eq:nepvdef} contains an orthogonality condition.

\section*{Acknowledgement}
The work of the second author was partially carried out during a research visit supported by EPFL/Univ. Geneva.
\bibliographystyle{siamplain}
\bibliography{eliasbib}

\begin{thebibliography}{10}

\bibitem{AMGHIBECH:2003}
{\sc S.~Amghibech}, {\em Eigenvalues of the discrete p-laplacian for graphs.},
  Ars Comb., 67 (2003), pp.~283--302.

\bibitem{Bai:2020:OPTIMALSCF}
{\sc Z.~Bai, R.-C. Li, and D.~Lu}, {\em Optimal convergence rate of
  self-consistent field iteration for solving eigenvector-dependent nonlinear
  eigenvalue problems}, tech. report, 2020,
  \url{https://arxiv.org/abs/2009.09022}.

\bibitem{Bao:2004:BOSEEINSTEIN}
{\sc W.~Bao and Q.~Du}, {\em Computing the ground state solution of
  {Bose-Einstein} condensates by a normalized gradient flow}, SIAM J. Sci.
  Comput., 25 (2004), pp.~1674--1697,
  \url{https://doi.org/10.1137/S1064827503422956}.

\bibitem{boyd2018simplified}
{\sc Z.~M. Boyd, E.~Bae, X.-C. Tai, and A.~L. Bertozzi}, {\em Simplified energy
  landscape for modularity using total variation}, SIAM J. Appl. Math., 78
  (2018), pp.~2439--2464.

\bibitem{BRIDLE:SOFTMAX:1989}
{\sc J.~Bridle}, {\em Training stochastic model recognition algorithms as
  networks can lead to maximum mutual information estimation of parameters}, in
  Adv. Neural Inf. Process. Syst., D.~Touretzky, ed., vol.~2, Morgan-Kaufmann,
  1990,
  \url{https://proceedings.neurips.cc/paper/1989/file/0336dcbab05b9d5ad24f4333c7658a0e-Paper.pdf}.

\bibitem{Hein:2009:PLAPLACIAN}
{\sc T.~B\"{u}hler and M.~Hein}, {\em Spectral clustering based on the graph
  p-{Laplacian}}, in Proceedings of the 26th International Conference on
  Machine Learning, 2009, pp.~81--88.

\bibitem{CAI:2018}
{\sc Y.~Cai, L.-H. Zhang, Z.~Bai, and R.-C. Li}, {\em On an
  eigenvector-dependent nonlinear eigenvalue problem}, SIAM J. Matrix Anal.
  Appl., 39 (2018), pp.~1360--1382, \url{https://doi.org/10.1137/17M115935X}.

\bibitem{Cances:2000:SCF}
{\sc E.~Canc\`es and C.~L. Bris}, {\em On the convergence of {SCF} algorithms
  for the {Hartree-Fock} equations}, M2AN, Math. Model. Numer. Anal., 34
  (2000), p.~749–774, \url{https://doi.org/10.1051/m2an:2000102}.

\bibitem{CHEN:1996}
{\sc C.~Chen and O.~L. Mangasarian}, {\em A class of smoothing functions for
  nonlinear and mixed complementarity problems}, Comp. Optim. and Appl., 5
  (1996), pp.~97--138, \url{https://doi.org/10.1007/BF00249052}.

\bibitem{SOFTPLUS:2000}
{\sc C.~Dugas, Y.~Bengio, F.~B\'{e}lisle, C.~Nadeau, and R.~Garcia}, {\em
  Incorporating second-order functional knowledge for better option pricing},
  in Adv. Neural Inf. Process. Syst., T.~Leen, T.~Dietterich, and V.~Tresp,
  eds., vol.~13, MIT Press, 2001,
  \url{https://proceedings.neurips.cc/paper/2000/file/44968aece94f667e4095002d140b5896-Paper.pdf}.

\bibitem{fasino2014algebraic}
{\sc D.~Fasino and F.~Tudisco}, {\em An algebraic analysis of the graph
  modularity}, SIAM J. Matrix Anal. Appl., 35 (2014), pp.~997--1018.

\bibitem{fortunato2016community}
{\sc S.~Fortunato and D.~Hric}, {\em Community detection in networks: A user
  guide}, Phys. Rep., 659 (2016), pp.~1--44.

\bibitem{HAGEN:1991}
{\sc L.~W. Hagen and A.~B. Kahng}, {\em Fast spectral methods for ratio cut
  partitioning and clustering}, 1991 IEEE International Conference on
  Computer-Aided Design Digest of Technical Papers,  (1991), pp.~10--13.

\bibitem{Hein:2010:IPM}
{\sc M.~Hein and T.~B\"{u}hler}, {\em An inverse power method for nonlinear
  eigenproblems with applications in 1-spectral clustering and sparse {PCA}},
  in Adv. Neural Inf. Process. Syst. 23, 2010, pp.~847--855.

\bibitem{hein2011beyond}
{\sc M.~Hein and S.~Setzer}, {\em Beyond spectral clustering-tight relaxations
  of balanced graph cuts.}, in Adv. Neural Inf. Process. Syst., 2011,
  pp.~2366--2374.

\bibitem{HOLLAND:1983}
{\sc P.~W. Holland, K.~B. Laskey, and S.~Leinhardt}, {\em Stochastic
  blockmodels: First steps}, Soc. Netw., 5 (1983), pp.~109--137,
  \url{https://doi.org/https://doi.org/10.1016/0378-8733(83)90021-7}.

\bibitem{IMPLICIT:2021}
{\sc E.~Jarlebring and P.~Upadhyaya}, {\em Implicit algorithms for eigenvector
  nonlinearities}, Numer. Algorithms,  (2021),
  \url{https://doi.org/10.1007/s11075-021-01189-4}.

\bibitem{jost2014nonlinear}
{\sc L.~Jost, S.~Setzer, and M.~Hein}, {\em Nonlinear eigenproblems in data
  analysis: {B}alanced graph cuts and the ratiodca-prox}, in Extraction of
  quantifiable information from complex systems, Springer, 2014, pp.~263--279.

\bibitem{keller2016general}
{\sc M.~Keller and D.~Mugnolo}, {\em General cheeger inequalities for
  $p$-laplacians on graphs}, Nonlinear. Anal. Theory. Methods. Appl., 147
  (2016), pp.~80--95.

\bibitem{kloster2014heat}
{\sc K.~Kloster and D.~F. Gleich}, {\em Heat kernel based community detection},
  in Proc. ACM SIGKDD Int. Conf. Knowl. Discov. Data Min., 2014,
  pp.~1386--1395.

\bibitem{krzakala2013spectral}
{\sc F.~Krzakala, C.~Moore, E.~Mossel, J.~Neeman, A.~Sly, L.~Zdeborov{\'a}, and
  P.~Zhang}, {\em Spectral redemption in clustering sparse networks}, Proc. of
  Natl. Acad. of Sci., 110 (2013), pp.~20935--20940.

\bibitem{Luo:2010}
{\sc D.~Luo, H.~Huang, and C.~Ding}, {\em On the eigenvectors of the
  $p$-laplacian}, Mach. Learn., 81 (2010), pp.~37--51,
  \url{https://doi.org/10.1007/s10994-010-5201-z}.

\bibitem{vonLuxburg:2007:SPECLUST}
{\sc U.~V. Luxburg}, {\em A tutorial on spectral clustering}, Stat. Comput.,
  (2007), pp.~395--417, \url{https://doi.org/10.1007/s11222-007-9033-z}.

\bibitem{NINA:2007:SOCIAL}
{\sc N.~Mishra, R.~Schreiber, I.~Stanton, and R.~E. Tarjan}, {\em Clustering
  social networks}, in Algorithms and Models for the Web-Graph, Berlin,
  Heidelberg, 2007, Springer Berlin Heidelberg, pp.~56--67.

\bibitem{newman2006finding}
{\sc M.~E. Newman}, {\em Finding community structure in networks using the
  eigenvectors of matrices}, Phys. Rev. E, 74 (2006), p.~036104.

\bibitem{ng2002spectral}
{\sc A.~Y. Ng, M.~I. Jordan, and Y.~Weiss}, {\em On spectral clustering:
  Analysis and an algorithm}, in Adv. Neural Inf. Process. Syst., 2002,
  pp.~849--856.

\bibitem{PASADAKIS:2020}
{\sc D.~Pasadakis, C.~L. Alappat, O.~Schenk, and G.~Wellein}, {\em $k$-way
  $p$-spectral clustering on grassmann manifolds}, tech. report, 2020,
  \url{https://arxiv.org/abs/2008.13210}.

\bibitem{Pulay:1980:CONVERGENCE}
{\sc P.~Pulay}, {\em Convergence acceleration of iterative sequences - the case
  of {SCF} iteration}, Chem. Phys. Lett, 73 (1980), pp.~393--398,
  \url{https://doi.org/10.1016/0009-2614(80)80396-4}.

\bibitem{REMM:2001}
{\sc M.~Remm, C.~E. Storm, and E.~L. Sonnhammer}, {\em Automatic clustering of
  orthologs and in-paralogs from pairwise species comparisons}, J. Mol. Bol.,
  314 (2001), pp.~1041--1052, \url{https://doi.org/10.1006/jmbi.2000.5197}.

\bibitem{Saad:2010:ELECSTRUCT}
{\sc Y.~Saad, J.~T. Chelikowsky, and S.~M. Shontz}, {\em Numerical methods for
  electronic structure calculations of materials}, SIAM Rev., 52 (2010),
  pp.~3--54, \url{https://doi.org/10.1137/060651653}.

\bibitem{SAUNDERS:LEVEL:1993}
{\sc V.~R. Saunders and I.~H. Hillier}, {\em A “level–shifting” method
  for converging closed shell hartree–fock wave functions}, Int. J. Quant.
  Chem., 7 (1973), pp.~699--705,
  \url{https://doi.org/https://doi.org/10.1002/qua.560070407}.

\bibitem{SCHMIDT:2007}
{\sc M.~Schmidt, G.~Fung, and R.~Rosales}, {\em Fast optimization methods for
  l1 regularization: A comparative study and two new approaches}, in Machine
  Learning: ECML 2007, 2007, pp.~286--297.

\bibitem{Selvan2017}
{\sc A.~N. Selvan, L.~M. Cole, L.~Spackman, S.~Naylor, and C.~Wright}, {\em
  Hierarchical Cluster Analysis to Aid Diagnostic Image Data Visualization of
  MS and Other Medical Imaging Modalities}, Springer New York, New York, NY,
  2017, pp.~95--123, \url{https://doi.org/10.1007/978-1-4939-7051-3_10}.

\bibitem{SHI:2000}
{\sc J.~Shi and J.~Malik}, {\em Normalized cuts and image segmentation}, EEE
  Trans. Pattern Anal. Mach. Intell., 22 (2000), pp.~888--905,
  \url{https://doi.org/10.1109/34.868688}.

\bibitem{szlam2010total}
{\sc A.~Szlam and X.~Bresson}, {\em Total variation, cheeger cuts}, in Intl.
  Conf. on Mach. Learn., 2010.

\bibitem{tudisco2018nodal}
{\sc F.~Tudisco and M.~Hein}, {\em A nodal domain theorem and a higher-order
  cheeger inequality for the graph $p$-laplacian}, J. Spectr. Theory., 8
  (2018), pp.~883--908.

\bibitem{tudisco2018community}
{\sc F.~Tudisco, P.~Mercado, and M.~Hein}, {\em Community detection in networks
  via nonlinear modularity eigenvectors}, SIAM J. Appl. Math., 78 (2018),
  pp.~2393--2419.

\bibitem{DENSITY:2021}
{\sc P.~Upadhyaya, E.~Jarlebring, and E.~H. Rubensson}, {\em A density matrix
  approach to the convergence of the self-consistent field iteration}, Numer.
  Algebra, Control. Optim., 11 (2021), pp.~99--115.

\bibitem{VORONIN:2015}
{\sc S.~Voronin, G.~Ozkaya, and D.~Yoshida}, {\em Convolution based smooth
  approximations to the absolute value function with application to non-smooth
  regularization}, tech. report, 2015, \url{https://arxiv.org/abs/1408.6795}.

\bibitem{WALKER:ANDERSON:2011}
{\sc H.~F. Walker and P.~Ni}, {\em Anderson acceleration for fixed-point
  iterations}, SIAM J. Numer. Anal., 49 (2011), pp.~1715--1735,
  \url{https://doi.org/10.1137/10078356X}.

\bibitem{XU:2015}
{\sc D.~Xu and Y.~Tian}, {\em A comprehensive survey of clustering algorithms},
  Ann. Data. Sci., 2 (2015), pp.~165--193,
  \url{https://doi.org/10.1007/s40745-015-0040-1}.

\bibitem{Yang:2009:SCF}
{\sc C.~Yang, W.~Gao, and J.~C. Meza}, {\em On the convergence of the
  self-consistent field iteration for a class of nonlinear eigenvalue
  problems}, SIAM J. Matrix Anal. Appl., 30 (2009), pp.~1773--1788,
  \url{https://doi.org/10.1137/080716293}.

\bibitem{YU:2020:PECOS}
{\sc H.-F. Yu, K.~Zhong, and I.~S. Dhillon}, {\em Pecos: Prediction for
  enormous and correlated output spaces}, tech. report, 2020,
  \url{https://arxiv.org/abs/2010.05878}.

\bibitem{ZERNER:DAMPING:1979}
{\sc M.~C. Zerner and M.~Hehenberger}, {\em A dynamical damping scheme for
  converging molecular scf calculations}, Chem. Phys. Lett., 62 (1979),
  pp.~550--554,
  \url{https://doi.org/https://doi.org/10.1016/0009-2614(79)80761-7}.

\end{thebibliography}

\end{document}